\documentclass{article}
\usepackage{amssymb,amsmath, euscript}

\input{epsf}

\newcounter{sec}

\def\SS{\smallskip}

\newcounter{punct}[sec]

\def\punct{\refstepcounter{punct}{\arabic{sec}.\arabic{punct}.  }}

 \newcounter{Apunct}
 
\def\Apunct{\refstepcounter{Apunct}{A.\arabic{Apunct}.  }}

\def\COUNTERS{\addtocounter{sec}{1}
              \setcounter{punct}{0}
          \setcounter{equation}{0}
          \setcounter{theorem}{0}
               }

\newtheorem{theorem}{Theorem}[sec]
\newtheorem{proposition}[theorem]{Proposition}

\newtheorem{lemma}[theorem]{Lemma}

\newtheorem{observation}[theorem]{Observation}

\newcounter{addendums}

\newtheorem{Atheorem}{Theorem}[addendums]

\begin{document}

\def\ov{\overline}

\def\wt{\widetilde}

\newcommand{\rk}{\mathop {\mathrm {rk}}\nolimits}
\newcommand{\Aut}{\mathop {\mathrm {Aut}}\nolimits}
\newcommand{\Out}{\mathop {\mathrm {Out}}\nolimits}
\renewcommand{\Re}{\mathop {\mathrm {Re}}\nolimits}

\def\Br{\mathrm {Br}}

\def\SL{\mathrm {SL}}
\def\SU{\mathrm {SU}}
\def\GL{\mathrm  {GL}}
\def\U{\mathrm  U}
\def\OO{\mathrm  O}
\def\Sp{\mathrm  {Sp}}
\def\SO{\mathrm  {SO}}
\def\SOS{\mathrm {SO}^*}
\def\Diff{\mathrm{Diff}}
\def\Vect{\mathfrak{Vect}}

\def\PGL{\mathrm  {PGL}}
\def\PU{\mathrm {PU}}

\def\PSL{\mathrm  {PSL}}

\def\Symp{\mathrm{Symp}}

\def\End{\mathrm{End}}
\def\Mor{\mathrm{Mor}}
\def\Aut{\mathrm{Aut}}

\def\PB{\mathrm{PB}}

\def\cA{\mathcal A}
\def\cB{\mathcal B}
\def\cC{\mathcal C}
\def\cD{\mathcal D}
\def\cE{\mathcal E}
\def\cF{\mathcal F}
\def\cG{\mathcal G}
\def\cH{\mathcal H}
\def\cJ{\mathcal J}
\def\cI{\mathcal I}
\def\cK{\mathcal K}
\def\cL{\mathcal L}
\def\cM{\mathcal M}
\def\cN{\mathcal N}
\def\cO{\mathcal O}
\def\cP{\mathcal P}
\def\cQ{\mathcal Q}
\def\cR{\mathcal R}
\def\cS{\mathcal S}
\def\cT{\mathcal T}
\def\cU{\mathcal U}
\def\cV{\mathcal V}
\def\cW{\mathcal W}
\def\cX{\mathcal X}
\def\cY{\mathcal Y}
\def\cZ{\mathcal Z}


\def\0{{\ov 0}}
\def\1{{\ov 1}}


\def\frA{\mathfrak A}
\def\frB{\mathfrak B}
\def\frC{\mathfrak C}
\def\frD{\mathfrak D}
\def\frE{\mathfrak E}
\def\frF{\mathfrak F}
\def\frG{\mathfrak G}
\def\frH{\mathfrak H}
\def\frJ{\mathfrak J}
\def\frK{\mathfrak K}
\def\frL{\mathfrak L}
\def\frM{\mathfrak M}
\def\frN{\mathfrak N}
\def\frO{\mathfrak O}
\def\frP{\mathfrak P}
\def\frQ{\mathfrak Q}
\def\frR{\mathfrak R}
\def\frS{\mathfrak S}
\def\frT{\mathfrak T}
\def\frU{\mathfrak U}
\def\frV{\mathfrak V}
\def\frW{\mathfrak W}
\def\frX{\mathfrak X}
\def\frY{\mathfrak Y}
\def\frZ{\mathfrak Z}

\def\fra{\mathfrak a}
\def\frb{\mathfrak b}
\def\frc{\mathfrak c}
\def\frd{\mathfrak d}
\def\fre{\mathfrak e}
\def\frf{\mathfrak f}
\def\frg{\mathfrak g}
\def\frh{\mathfrak h}
\def\fri{\mathfrak i}
\def\frj{\mathfrak j}
\def\frk{\mathfrak k}
\def\frl{\mathfrak l}
\def\frm{\mathfrak m}
\def\frn{\mathfrak n}
\def\fro{\mathfrak o}
\def\frp{\mathfrak p}
\def\frq{\mathfrak q}
\def\frr{\mathfrak r}
\def\frs{\mathfrak s}
\def\frt{\mathfrak t}
\def\fru{\mathfrak u}
\def\frv{\mathfrak v}
\def\frw{\mathfrak w}
\def\frx{\mathfrak x}
\def\fry{\mathfrak y}
\def\frz{\mathfrak z}

\def\frsp{\mathfrak{sp}}


\def\bfa{\mathbf a}
\def\bfb{\mathbf b}
\def\bfc{\mathbf c}
\def\bfd{\mathbf d}
\def\bfe{\mathbf e}
\def\bff{\mathbf f}
\def\bfg{\mathbf g}
\def\bfh{\mathbf h}
\def\bfi{\mathbf i}
\def\bfj{\mathbf j}
\def\bfk{\mathbf k}
\def\bfl{\mathbf l}
\def\bfm{\mathbf m}
\def\bfn{\mathbf n}
\def\bfo{\mathbf o}
\def\bfp{\mathbf p}
\def\bfq{\mathbf q}
\def\bfr{\mathbf r}
\def\bfs{\mathbf s}
\def\bft{\mathbf t}
\def\bfu{\mathbf u}
\def\bfv{\mathbf v}
\def\bfw{\mathbf w}
\def\bfx{\mathbf x}
\def\bfy{\mathbf y}
\def\bfz{\mathbf z}

\def\bfA{\mathbf A}
\def\bfB{\mathbf B}
\def\bfC{\mathbf C}
\def\bfD{\mathbf D}
\def\bfE{\mathbf E}
\def\bfF{\mathbf F}
\def\bfG{\mathbf G}
\def\bfH{\mathbf H}
\def\bfI{\mathbf I}
\def\bfJ{\mathbf J}
\def\bfK{\mathbf K}
\def\bfL{\mathbf L}
\def\bfM{\mathbf M}
\def\bfN{\mathbf N}
\def\bfO{\mathbf O}
\def\bfP{\mathbf P}
\def\bfQ{\mathbf Q}
\def\bfR{\mathbf R}
\def\bfS{\mathbf S}
\def\bfT{\mathbf T}
\def\bfU{\mathbf U}
\def\bfV{\mathbf V}
\def\bfW{\mathbf W}
\def\bfX{\mathbf X}
\def\bfY{\mathbf Y}
\def\bfZ{\mathbf Z}

\def\bfw{\mathbf w}

\def\R {{\mathbb R }}
 \def\C {{\mathbb C }}
  \def\Z{{\mathbb Z}}
  \def\H{{\mathbb H}}
\def\K{{\mathbb K}}
\def\N{{\mathbb N}}
\def\Q{{\mathbb Q}}
\def\A{{\mathbb A}}

\def\T{\mathbb T}
\def\P{\mathbb P}

\def\G{\mathbb G}

\def\bbA{\mathbb A}
\def\bbB{\mathbb B}
\def\bbD{\mathbb D}
\def\bbE{\mathbb E}
\def\bbF{\mathbb F}
\def\bbG{\mathbb G}
\def\bbI{\mathbb I}
\def\bbJ{\mathbb J}
\def\bbL{\mathbb L}
\def\bbM{\mathbb M}
\def\bbN{\mathbb N}
\def\bbO{\mathbb O}
\def\bbP{\mathbb P}
\def\bbQ{\mathbb Q}
\def\bbS{\mathbb S}
\def\bbT{\mathbb T}
\def\bbU{\mathbb U}
\def\bbV{\mathbb V}
\def\bbW{\mathbb W}
\def\bbX{\mathbb X}
\def\bbY{\mathbb Y}

\def\kappa{\varkappa}
\def\epsilon{\varepsilon}
\def\phi{\varphi}
\def\le{\leqslant}
\def\ge{\geqslant}

         \begin{center}
\bf\Large

Infinite tri-symmetric group,\\
 multiplication of double cosets,
\\and
checker topological field theories

\bigskip
\sc\large{Yury Neretin}%
\footnote{Supported by grant   FWF, project P19064}
\end{center}

{\small We consider a product of three copies of infinite
symmetric group and its representations spherical
 with respect to the diagonal subgroup. We show that such
 representations generate functors from a certain category
 of simplicial two-dimensional surfaces to the category
  of Hilbert spaces
 and bounded linear operators.}

\section{Introduction}

\COUNTERS

{\bf\punct Infinite symmetric group.}
Denote by $S_\infty$ the group of all permutations
of the set $\N=\{1,2,3,\dots\}$. By $S_\infty(\alpha)$
 we denote the stabilizer
of  points $1,$ $2$, \dots, $\alpha\in \N$.
We also assume $S_\infty(0)=S_\infty$.

 A topology on the group $S_\infty$
is defined by the condition: subgroups $S_\infty(\alpha)$ are open
and
cosets $gS_\infty(\alpha)\subset S_\infty$ form a basis of topology.
 The group $S_\infty$ is a totally disconnected topological
group.

Classification of irreducible
unitary representations of $S_\infty$ is rather
simple, all representations are induced from trivial representations
of Young subgroups
$S_{m_1}\times\dots\times S_{m_k}\times S_{\infty-\sum m_j}$,
see
Lieberman   \cite{Lie}, see also
Olshanski \cite{Ols-lieb}, and exposition in \cite{Ner-book},
VIII,1-2.

\smallskip

By $S_\infty^{fin}\subset S_\infty$
we denote the group of finite permutations
of $\N$, this group is an inductive limit
of finite symmetric groups,
$$
S_\infty^{fin}=\cup_{n=1}^\infty S_n
.
$$

\smallskip

The group $S_\infty^{fin}$ is a wild (not type I, see e.g.,
\cite{Dix}) discrete group,
and description of its representations seems to be a non-reasonable
problem.

\smallskip


{\bf\punct $n$-symmetric groups.} Consider the product
$$
 G^{[n]}:=S_\infty\times \dots \times S_\infty
$$
 of $n$ copies of $S_\infty$. We write  elements $g\in  G^{[n]}$ as
 collections $\bfg=(g_1,\dots, g_n)$, where $g_j\in S_\infty$.
Denote by $K$ the diagonal subgroup in $G^{[n]}$,
 i.e., the subgroup
 consisting of collections $(g,g,\dots,g)$, where
 $g\in S_\infty$.

We define the {\it $n$-symmetric group} $\G=\G^{[n]}$
 as the subgroup of
$G^{[n]}$ consisting of collections $(g_1,\dots,g_n)$
such that
$$
g_i g_j^{-1}\in S_\infty^{fin}
\qquad\text{for all $i$, $j\le n$}
.
$$
In other words, $\G^{[n]}$ consists of all collections
$$
(g h_1, gh_2,\dots,g h_n) \qquad\text{such that $g\in S_\infty$,
 $h_j\in S_\infty^{fin}$}
.
$$

Denote by $K(\alpha)$ the image of $S_\infty(\alpha)$ under
 the diagonal embedding $K\mapsto \G^{[n]}$.
 We define the topology on $\G^{[n]}$ from the condition:
 the subgroups $K(\alpha)$ are open.
In other words, the topology of $K\simeq S_\infty$
is the same as above,
 the quotient-space $\G^{[n]}/K$ is countable and equipped with
 the discrete topology.

\smallskip


{\bf\punct Bisymmetric group.} The existing representation theory
of  infinite symmetric groups is mainly the representation theory
of the bisymmetric group $\G^{[2]}$, see Thoma \cite{Tho},
Vershik, Kerov \cite{VK}, Olshanski \cite{Olsh-symm},
 Okounkov \cite{Oko}, Kerov, Olshanski, Vershik
 \cite{KOV}.  The situation was explained by
 Olshanski in \cite{Olsh-symm}. We refer the reader to this paper.

The group $\G^{[n]}$ is outside Olshanski's approach to
infinite-dimensional groups, based on imitation of symmetric pairs,
see Olshanski \cite{Ols-howe}, \cite{Olsh-symm}, and also
\cite{Ner-book}. However, $\G^{[n]}$ is a $(G,K)$-pair in the sense
of \cite{Ner-book}, VIII.5. My another standpoint was a strange
Nessonov's theorem, \cite{Ness1}, \cite{Ness2} about spherical
functions on certain infinite-dimensional groups
 with respect to certain small (non-symmetric) subgroups.

A discussion of bisymmetric groups is contained
in Addendum to this paper.


\smallskip

{\bf\punct Content of the paper. The tri-symmetric group.}
In Sections 2--4 we discuss the tri-symmetric group
$\G:=\G^{[3]}$.

\smallskip

 Section 2. We describe of double cosets%
 \footnote{Let $G$ be a group, $H$, $K$ its subgroups.
 A double coset  is a subset in $G$
  of the form $HgK$, where $g\in G$.
 The space $H\setminus G/K$ is the set, whose points
 are double cosets. Denote by $\cF(H\setminus G/K)$
 the space of functions on $H\setminus G/K$
  If a group $G$ is finite, then there is
 a well-defined convolution
$\cF(H_1\setminus G/H_2)\times \cF(H_2\setminus G/H_3)
\to \cF(H_1\setminus G/H_3)$. But a product of individual
double cosets is not well-defined.}.
 The space $H\setminus G/K$
$K(\alpha)\setminus\G/K(\beta)$.
 More precisely, we identify a double coset
with a triangulated compact two-dimensional surface equipped
with certain additional data ('checker-boards').

\smallskip

 Next, we construct a natural multiplication
$(\fra,\frb)\mapsto\fra\circ\frb$ of double cosets,
$$
K(\alpha)\setminus\G/K(\beta)
\,\times\, K(\beta)\setminus \G/K(\gamma)
\,\to\,
K(\alpha)\setminus \G/K(\beta)
$$
for each $\alpha$, $\beta$, $\gamma\in\Z_+$. Thus we get
a category $\bbS=\bbS^{[3]}$, whose objects are $0$, $1$, $2$,
 \dots
and morphisms $\beta\to\alpha$ are double cosets
$K(\alpha)\setminus \G/K(\beta)$.

Multiplications of double cosets is a usual phenomenon for
infinite-dimensional groups
(Ismagilov--Olshanski multiplicativity, see numerous
constructions
in \cite{Ner-book}, see also some additional examples
 in Russian
translation of \cite{Ner-book}, Addendum E, and in
\cite{Ner-poli}).

In this case, we get a construction similar to
so-called 'topological field theories'%
\footnote{Topological field theories
are  imitations of conformal field theories
on topological level. Usually, this term is used for functors
from categories of bordisms to the category of linear
spaces and linear operators. Bordisms of simplicial
complexes were considered, e.g., by Natanzon \cite{Nat}.}.

\smallskip

Section 3. Let $\rho$ be a unitary  representation
 of $\G$ in a Hilbert space
$H$. Denote by $H(\alpha)$ the space of $K(\alpha)$-fixed vectors.
Then each $\fra\in K_\alpha\setminus\G/ K_\beta$ determines
 a well-defined operator $\rho(\fra):H(\beta)\to H(\alpha)$
 and
 $$
\rho(\fra)\rho(\frb)=\rho(\fra\circ\frb)
\qquad\text{for
$\fra\in K(\alpha)\setminus\G/K(\beta)$,\,
$\frb\in K(\beta)\setminus \G/K(\gamma)$,}
 $$
i.e., we get a representation of the category $\bbS$.

\smallskip

Section 4. We construct a family of representations of $\bbG$,
write expressions for their spherical functions
and for matrix elements  of representations of the category
$\bbS$.

\smallskip


{\bf\punct $n$-symmetric groups.}
General $n$-symmetric groups
$\G^{[n]}$  are similar to the tri-symmetric
group. Generally,  checker-boards are tiled by
$n$-gons. Also, there is another (nonequivalent) construction
of a two-dimensional polygonal complex from a double coset.
We discuss this in Section 5.

\smallskip


{\bf\punct Abstract theorems.}
Sections 1--5 form the main part of the paper. In
supplementary Section 6, we obtain some
abstract theorems about unitary representations of
$n$-symmetric groups: uniqueness of a $K$-spherical vector,
rough separation of the set of representations,
etc.

\smallskip


{\bf Acknowledgements.} I am grateful to G.~I.~Olshanski
and N.~I.~Nessonov for discussions of this topic.

\section{Multiplication of checker-boards}

\COUNTERS

Here we consider the tri-symmetric group $\G:=\G^{[3]}$.

\smallskip

{\bf \punct Construction of simplicial complexes.}
We have 3 copies of the set $\N$, say {\it red, yellow, and
blue.} An element of $\G$ is  a
triple of permutations of $\N$, denote it by
$$
(g_{red}, \,\,g_{yellow},\,\, g_{blue}).$$

 We draw a  collection of disjoint oriented
 {\it black} triangles $A_j$, where $j$ ranges in $\N$,
and paint their  sides in red, yellow, and blue {\it anti-clockwise}.
We also draw a  collection of oriented
 {\it white}
triangles $B_j$ and  paint sides in red, yellow, and
 blue {\it clockwise}.

Next, we glue a simplicial complex from these triangles.
If $g_{red}$  sends $i$ to $j$, then we identify
the red side of the  black triangle $A_i$ with the  red side
of the white triangle $B_{j}$. We repeat the same
operation for $g_{yellow}$ and $g_{blue}$.
 In this way, we get a disjoint
countable union
of 2-dimensional compact  closed triangulated surfaces.

\begin{figure}
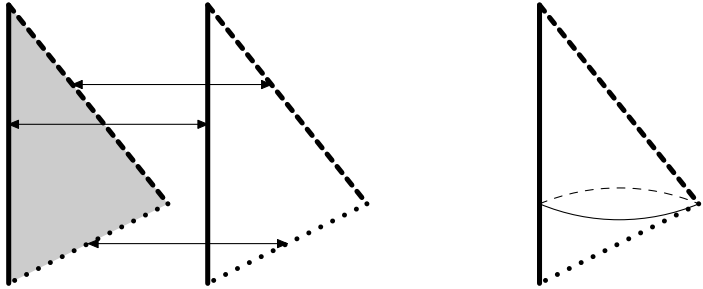

\epsfbox{tri.2} \qquad\qquad\qquad \epsfbox{tri.3}

\caption{Gluing of a cheburek}
\label{fig:cheburek}

\end{figure}

All components except finite number consist of two triangles
(black and white, glued along the corresponding sides).
We call such components {\it chebureks}%
\footnote{A street fast-food in some countries.},
 see Figure
\ref{fig:cheburek}.

\smallskip

{\sc Remark.} Elements of the subgroup $K\subset\G$
correspond to countable collections of disjoint chebureks.

\smallskip


\begin{figure}
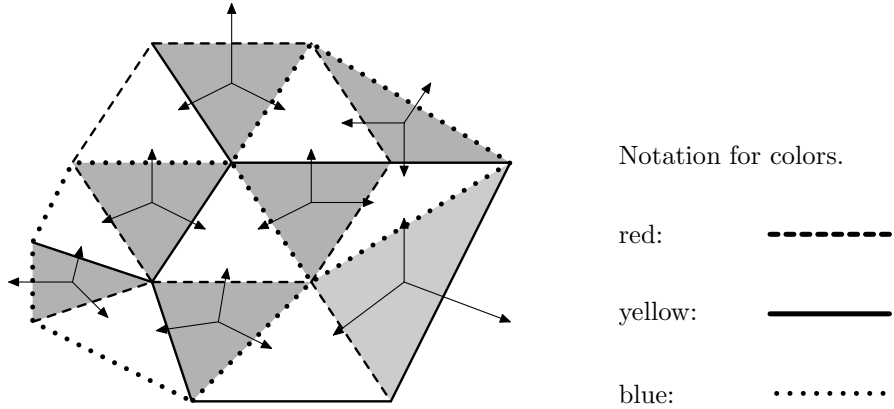


$$
\epsfbox{tri.1}\qquad \qquad \epsfbox{tri.9}
$$

\caption{A checker-board and 3 substitutions}
\label{fig:board}

\end{figure}


{\bf\punct Checker-boards.}
Now we define a {\it checker-board}
(see Figure \ref{fig:board}) as a countable disjoint union
of triangulated closed surfaces equipped  with following data:

\smallskip

a) Triangles are painted in black and white, neighbors
of black triangles  are white and vise versa.

\smallskip

b) Edges are painted in red,  yellow, and blue. The boundary
of any black (respectively white) triangles is composed of red,
 yellow, and blue edges situated
  anti-clockwise (resp., clockwise).

\smallskip

c) Black (respectively, white) triangles are enumerated
by natural numbers.

\smallskip

d) Almost all components are chebureks,
i.e., unions of two triangles.

\begin{proposition}
There is a canonical one-to-one correspondence between
$\G$ and the set of all checker-boards.
\end{proposition}

One-side construction was given above.
Conversely, take a checker-board. Let a red (yellow, blue) edge
be common for $i$-th black triangle and $j$-th white triangle.
Then $g_{red}$ sends $i$ to $j$, see Figure \ref{fig:board}.

\smallskip


{\bf\punct Vertices of checker-boards.}
There are 3 types of vertices according colors of edges at
a vertex: red--blue, red--yellow, yellow--blue.

\begin{proposition}
The set of red--blue vertices is in
a natural one-to-one correspondence
with the set of independent cycles of $g_{red}g^{-1}_{blue}$.
Moreover, the number of triangles meeting at a vertex
is the double length of the corresponding cycle.
\end{proposition}

This is obvious, see Figure \ref{fig:board}.


\smallskip

{\bf\punct Double cosets and checker-boards.}
We say that an {\it $(\alpha,\beta)$-board} is a compact
(generally, disconnected) oriented  triangulated surface
equipped
with the following data:

\SS

a) triangles are painted in white and black,
edges are painted in red, yellow, blue as
above,

\smallskip

b) we assign labels 1, 2, \dots, $\beta$
to some black triangles and 1, 2, \dots, $\alpha$ to some white
triangles.

\SS

c) There is no chebureks without labels.

\begin{proposition}
There is a canonical one-to-one correspondence between the set
$K(\alpha)\setminus \G/K(\beta)$ and the set of all
$(\alpha,\beta)$-boards.
\end{proposition}

Indeed, take a double coset, choose
its representative $\in \G$ and
compose the checker-board for $g$.
Next, we remove numbers
$>\beta$ from black triangles and numbers $>\alpha$ from
white triangles. Now almost all components are chebureks
without numbers (we call them {\it empty chebureks}),
 and we remove them.

Thus we identify double cosets and checker-boards.


\smallskip

{\bf\punct Multiplication of checker-boards.}
Consider two checker-boards,
$$
\fra\in K(\alpha)\setminus\G /K(\beta)
,\qquad
\frb\in K(\beta)\setminus\G/K(\gamma)
.
$$
For each $j=1$, $2$, \dots, $\beta$ we cut off
the black triangle of $\fra$ with label
$j$ and the white triangle of $\frb$ with label $j$
and identify edges of the triangles according colours
and orientations.
We obtain an $(\alpha,\gamma)$-board.

  Denote this operation by $\fra\circ\frb$.

\begin{figure}
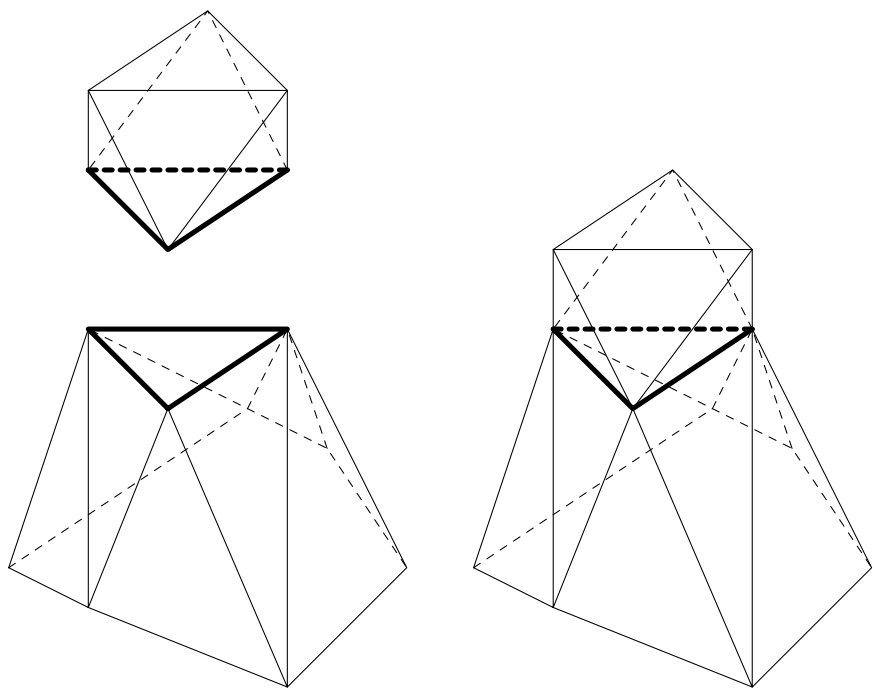

\epsfbox{tri.4} \qquad \epsfbox{tri.5}

\caption{Gluing of checker-boards. On the resulting surface,
the interior part of the fat triangle is removed.}
\end{figure}

Thus we get a category, whose objects are $0$, $1$, $2$, \dots,
and morphisms $\beta\to\alpha$
 are $(\alpha,\beta)$-boards or,
 equivalently, double cosets.

The real reason for this definition
is Theorem \ref{th:io} about representations.
 Now we wish to interpret the definition
in terms of the group $\G$.


\SS

{\bf\punct Weak zeros.%
\label{ss:weak}}
Let us write elements of $K\simeq S_\infty$ as 0-1 matrices.
We say that a sequence $h_j$ in $K(\alpha)$ is a {\it weak zero}
if $h_j$ converges element-wise to the
$(\alpha+\infty)\times(\alpha+\infty)$-matrix
$\begin{pmatrix} 1&0\\0&0 \end{pmatrix}$.

\smallskip


{\bf\punct  Another definition
of the multiplication of double cosets.}

\begin{theorem}
\label{th:ugolok}
 Given
$$\fra\in K(\alpha)\setminus \G/K(\beta),
\quad
\frb \in K(\beta)\setminus \G/K(\gamma)
,$$
choose $\bfp\in\fra$, $\bfq\in\frb$, and a weak zero $h_j$ in
$K(\beta)$.
Then
$$
\frc_j:=K(\alpha)\bfp h_j\bfq K(\gamma)
$$
equals  $\frc=\fra\circ\frb$  for sufficiently large $j$.
\end{theorem}

{\sc Proof.} First, consider two elements
$\bfp$, $\bfq\in\G$ and the corresponding checker-boards
$\frP$, $\frQ$.
 Then the product $\bfp\bfq$
in the group $\G$ corresponds to the following
operation:

\smallskip

 --- we remove black triangles from
$\frP$  and white triangles from $\frQ$;

\smallskip

--- we glue perimeters of corresponding 'holes'
in $\frP$ and $\frQ$ (according colors and orientations of edges).

\smallskip

Second, let $h\in K\subset\G$. Then the product
$\bfp h\bfq$ corresponds to the following operation:
for each $j$ we glue $j$-th white triangle
of $\frQ$ with $h\cdot j$-th black triangle of $\frP$.

Let $\bfp\in\G$.
Let draw the corresponding checker-board $\frp$.
Let us remove all chebureks, whose black side has label
$>\beta$ and white side has label $>\alpha$, denote by
$\cP$ the corresponding reduced checker board.
We say that the {\it right} (resp. {\it left})
 {\it $(\alpha,\beta)$-support} of $\bfp$ is the set
 of all labels on black (resp. white) triangles
 of   $\cP$.

  Equivalently, $j$ is not in the right support
 iff
 $$
 j>\beta \qquad\text{and}\qquad
p_{red} j=p_{yellow} j=p_{blue}>\beta
 .$$

Now, let $h_j$ be a weak zero in $K(\beta)$.
For sufficiently large $j$, the permutation
  $h_j\in K(\beta)$ sends the left support of $\bfq$
  to the complement of the right support of $\bfp$.

 Now let $\bfp$, $\bfq\in\G$, let $\frP$, $\frQ$
 be the corresponding checker-boards.
  We evaluate $\bfp h_j \bfq$
 by the rule described above. Then all black triangles
 of $\frQ$ with labels $>\beta$ are glued with chebureks,
 whose white labels are big. A gluing of a cheburek
 does not change a surface. Big label will be forgotten.

  The same holds for white triangles
 of $\frP$ with labels $>\beta$. Thus, we get the operation
 of gluing of $(\alpha,\beta)$-board and $(\beta,\gamma)$-board
  described above.
  \hfill $\square$

\smallskip


{\bf \punct Involution in the category $\bbS$.}
The map $\bfp\mapsto \bfp^{-1}$ induces a map
$$
K(\alpha)\setminus\G/K(\beta)\,
\to\,
K(\beta)\setminus \G/K(\alpha)
.$$
For an $(\alpha,\beta)$-board $\fra$ we make
a $(\beta,\alpha)$-board $\fra^\square$ by the follow rule:

\smallskip

a) we change 'black' and 'white';

\smallskip

b) we change the orientation of the surface.

\smallskip

Evidently,
$$
(\fra\circ \frb)^\square=\frb^\square\circ\fra^\square
.
$$


\section{Representations. Abstract theorem}

\COUNTERS

{\bf\punct Operators $\rho(\fra)$.}
Let $\rho$ be a unitary representation of $\G=\G^{[3]}$
in a Hilbert space $H$.
Denote by $H(\alpha)$ the space of $K(\alpha)$-fixed vectors.
Denote by $P(\alpha)$ the operator of orthogonal projection
 to the space
$H(\alpha)$.

Fix $\alpha$, $\beta$.
For $\bfp\in \G$, consider the operator
$$
H(\beta)\to H(\alpha)
$$
given by
$$
\ov \rho(g):=
P(\alpha) \rho(\bfp)=P(\alpha) \rho(\bfp) P(\beta)
.
$$
For $r_1\in K(\alpha)$, $r_2\in K(\beta)$ we have
\begin{equation}
\label{eq:biinvariance}
\ov\rho(r_1 g r_2)=\ov\rho(g)
.
\end{equation}
Thus $\ov \rho$ is a function on double cosets
$K(\alpha)\setminus \G/K(\beta)$.

\begin{theorem}
\label{th:io}
For any unitary representation
$\rho$ of $\G$ for each $\alpha$, $\beta$,
$\gamma\in \Z_+$, for each
\begin{equation}
\label{eq:dd}
\fra\in K(\alpha)\setminus \G/K(\beta),\quad
\frb \in K(\beta)\setminus \G/K(\gamma)
,
\end{equation}
we have
$$
\ov \rho(\fra)\ov \rho(\frb)=\ov \rho(\fra\circ \frb)
.$$
\end{theorem}

Thus we get a representation of the category $\bbS$
of checker-boards. Precisely, for each $\alpha$ we assign
the Hilbert space $H(\alpha)$ and for each $(\alpha,\beta)$-board
$\fra$ we assign the operator $\rho(\fra):H(\beta)\to H(\alpha)$.

Evidently, $\ov\rho$ is a {\it $*$-representation}, i.e.,
$$
\ov\rho(\fra^\square)=\rho(\fra)^*.
$$


{\bf\punct Lemma on representations of $S_\infty$.}
For representations of $S_\infty$, see
\cite{Lie}, \cite{Ols-new}, \cite{Olsh-symm},
\cite{Ner-book}. For proof of the following Proposition
\ref{l:weak},
see, e.g., \cite{Ner-book},  Corollary 8.1.5.

For  a unitary representation of $S_\infty$
in a Hilbert space $H$,
we define subspaces $H(\alpha)\subset H$ as above.

\begin{proposition}
\label{l:weak}
Let $h_j$ be a weak zero%
 \footnote{See Subsection \ref{ss:weak}}
 in $S_\infty(\alpha)$.
Then $\rho(h_j)$ converges to $P(\alpha)$ in the weak
operator topology.
\end{proposition}

 \smallskip


{\bf\punct Proof of Theorem \ref{th:io}.} Let $\fra$, $\frb$
be double cosets as above (\ref{eq:dd}),
 let $\bfp$, $\bfq\in\G$ be their
representatives. Let $u_j$ be a weak zero in $K(\beta)$.
Then
\begin{multline*}
\ov\rho(\fra)\ov\rho(\frb)=
P(\alpha)\rho(\bfp)P(\beta)\rho(\bfq)P(\gamma)
=\\=
\lim_{j\to\infty}P(\alpha)\rho(\bfp)\rho(u_j)\rho(\bfq)P(\gamma)
=
P(\alpha)\lim_{j\to\infty}\rho(\bfp u_j \bfq)P(\gamma)
\end{multline*}
where $\lim$ denotes the weak limit.

For sufficiently large  $j$,
a product $\bfp u_j \bfq$ is contained in the double coset
$\fra\circ\frb$. Therefore, we get $\ov\rho(\fra\circ\frb)$.
\hfill$\square$

\section{Representations. Constructions}

\COUNTERS

In this section we  construct some representations
of the tri-symmetric group and of the category $\bbS$.


\smallskip

{\bf\punct Tensor product construction.} Here we
imitate the construction of \cite{VK}, \cite{Olsh-symm}.
Let $U$, $V$, $W$ be (finite dimensional
or infinite dimensional) Hilbert spaces.
 Let $u_i$, $v_j$, $w_k$ be their
orthonormal  bases.
Consider the tensor product
$$
H:=
U\otimes V\otimes W
.
$$
Fix a unit vector $h\in H$,
$$
h=\sum_{i,j,k} h_{ijk}\,\,u_i\otimes v_j\otimes w_k
\in U\otimes V\otimes W, \qquad \|h\|=1
.
$$
Consider the  tensor product
$$
\bfH:=(H,h)\otimes (H,h)\otimes\dots=
(U\otimes V\otimes W,\,h)\,\otimes (U\otimes V\otimes W,\,h)\,
\otimes\dots
$$
of infinite number of copies of Hilbert spaces $H$ with
the the distinguished vector $h$ (see the von Neumann
definition of infinite tensor products, \cite{vN}).
Denote
$$
\xi:=h\otimes h\otimes\dots\in \bfH.
$$

Now we wish to define a certain representation
$\nu_h$ of the group $\G$ in the space $\bfH$.
Let $\bfp=(p_1,p_2,p_3)\in \G$. Then
$p_1$ acts as a permutation of factors $U$, $p_2$ as permutation
of factors $V$, and $p_3$ as a permutation of factors $W$.

Let formulate this more carefully.
The operators of permutations are well defined

\smallskip

--- for
$p_1$, $p_2$, $p_3\in S_\infty^{fin}$;

\smallskip

--- for elements $(r,r,r)\in K$, they
 act by simultaneous permutations of
factors $(U\otimes V\otimes W)$ and preserve the distinguished
vector $\xi$.

\smallskip

These operators generate a unitary representation  of $\G$.

\smallskip

{\sc Remark.} The group $S_\infty\times S_\infty\times S_\infty$
does not act in this tensor product. \hfill $\square$

\smallskip

{\sc Remark.} Let $A$, $B$, $C$ be unitary operators in
$U$, $V$, $W$ respectively. Then
$$
\nu_h\simeq \nu_{(A\otimes B\otimes C)h}
.
$$
However, we can not simplify an expression for $h$ by such
transformations.
\hfill$\square$


\smallskip

{\bf\punct Spherical functions.}
Our representations are $K$-spherical in the following sense:

\begin{proposition}
\label{pr:spherical}
 The vector $\xi$ is a unique $K$-fixed
vector in $\bfH$.
\end{proposition}

{\sc Proof.}
Indeed, choose an orthonormal basis $e_1$, $e_2$, \dots
$\in U\otimes V\otimes W$ with $e_1=h$.
Expand a $K$-fixed vector $\eta\in \bfH$ in Fourier series ,
$$
\eta=\sum_{i_1,i_2,\dots;\,\, \text{where $i_N=1$ for large $N$}}
c_{i_1 i_2\dots}\,\,
e_{i_1}\otimes e_{i_2}\otimes\dots
$$
Coefficients $c$ must be invariant with respect to all permutations
of indices. Let some $i_k\ne 1$ and $c_{i_1 i_2\dots}=\epsilon\ne0$.
Then there is a countable number of pairwise different permutations
of the sequence $i_1$, $i_2$, \dots and therefore a countable number
of Fourier coefficients $=\epsilon$. Therefore $\|\eta\|=\infty$.
\hfill $\square$

\smallskip

Next, we write
the spherical function, i.e.,
$$
\Phi_h(\bfg)=\langle
 \ov\nu_h(\bfg)  \xi, \xi
 \rangle.
$$
By (\ref{eq:biinvariance},
$\Phi_h$ is a function on $K(0)\setminus \G/K(0)$, i.e.,
on the set of non-labeled checker-boards.

Let us write natural numbers $i\le\dim U$ on red edges
of the simplicial complex,
 $j\le\dim V$ on yellow edges, and $k\le \dim W$
on blue edges. For each triangle $T$,
 we get  3 numbers on its sides, say
$i(T)$, $j(T)$, $k(T)$.

\begin{proposition}
\begin{multline}
\Phi_h(\fra)
=\\=
\sum\limits_{\begin{matrix}\scriptstyle\text{all arrangements of}
\\\scriptstyle
\text{$i$, $j$, $k$ on edges}\end{matrix}}
\,\,
\prod\limits_{\text{black triangles $T$}}
h_{i(T)j(T)k(T)}
\prod\limits_{\text{white triangles $S$}}
\ov h_{i(S)j(S)k(S)}
\label{eq:spherical}
\end{multline}
\end{proposition}

{\sc Proof.} We choose a representative $\bfg\in\G$
of the coset $\fra$
such that $g_1$, $g_2$, $g_3$ are finite,
let they are contained in $S_N$. Hence we can work
in the space $H^{\otimes N}$,
\begin{multline}
\label{eq:long-1}
h^{\otimes N}=
\sum\limits_{ i_1,i_2,\dots;\,\,
 j_1,j_2,\dots;\,\,
 k_1,k_2,\dots}
 \\
\Bigl(\prod_{m\le N}
h_{i_m,j_m,k_m}\Bigr)
(u_{i_1}\otimes v_{j_1}\otimes w_{k_1})
\otimes \dots\otimes
(u_{i_N}\otimes v_{j_N}\otimes w_{k_N})
\end{multline}
Thus summands are enumerated by collections of numbers
written
on sides of black triangles.

 Next,
\begin{multline}
\label{eq:long-2}
\rho(g)
h^{\otimes N}
=
\sum\limits_{ i_1,i_2,\dots;\,\, j_1,j_2,\dots;\,\,
  k_1,k_2,\dots}
\\
\Bigl(\prod_{m\le N}
h_{i_m,j_m,k_m}\Bigr)
(u_{g_1(i_1)}\otimes v_{g_2(j_1)}\otimes w_{g_3(k_1)})
\otimes \dots\otimes
(u_{g_1(i_N)}\otimes v_{g_2(j_N)}\otimes w_{g_3(k_N)})
\end{multline}
We change the order of summation,
$$
i_m^{new}=g_1^{-1}(i_m),\quad j_m^{new}=g_2^{-1}(j_m),\quad
k_m^{new}=g_3^{-1}(k_m)
$$
and get
\begin{multline*}
\sum\limits_{1,i_2,\dots,\,\,
 j_1,j_2,\dots,\,\,
 k_1,k_2,\dots} \\
\Bigl(\prod_{m\le N}
h_{g_1^{-1}(i_m),g_2^{-1}(j_m),g_3^{-1}(k_m)}\Bigr)
(u_{i_1}\otimes v_{j_1}\otimes w_{k_1})
\otimes \dots\otimes
(u_{i_N}\otimes v_{j_N}\otimes w_{k_N})
\end{multline*}
Evaluating the inner product
of (\ref{eq:long-1}) and (\ref{eq:long-2}),
we come to (\ref{eq:spherical}).
\hfill$\square$

\smallskip

{\sc Remark.} a) The multiplication in
$K(0)\setminus \G /K(0)$ is the disjoint union.
Therefore, for $\fra$, $\frb\in K(0)\setminus \G /K(0)$
we have
$$
\Phi_h(\fra\circ \frb)=\Phi_h(\fra) \Phi_h(\frb)
.
$$
Factors corresponding to chebureks are 1.

\smallskip

b)The function $\fra\mapsto\Phi_h(\fra)\Phi_{h'}(\fra)$
has the form $\Phi_{\wt h}(\fra)$. Indeed,
having spaces
$$
H=
\bigl(U\otimes V\otimes W, \, h\bigr)
\qquad\text{and}\qquad
H'=
\bigl(U'\otimes V'\otimes W',\, h'\bigr)
,
$$
we take the space
$$
H\otimes H'=
\Bigl( (U\otimes U')\otimes (V\otimes V')\otimes (W\otimes W'),
\,h\otimes h'\Bigr)
.
$$

\smallskip


{\bf\punct Representations of the category $\bbS$.}
Now we wish to construct the corresponding representations
of the category of checker-boards.

\begin{proposition}
The spaces of $K(\alpha)$-fixed vectors are
$$
\bfH(\alpha)=\underbrace{H\otimes\dots\otimes H}_{\text{
$\alpha$ times}}
\otimes h\otimes h\otimes \dots
.$$
\end{proposition}

This a rephrasing of Proposition \ref{pr:spherical}.

\smallskip

In particular, $\bfH(\alpha)$ are finite dimensional iff
$U$, $V$, $W$ are finite dimensional.

\smallskip

Now let us write matrix elements for a checker-board
$\fra\in K(\alpha)\setminus \G/K(\beta)$.
Consider basis vectors
$$
\xi:=(u_{i_1}\otimes v_{j_1}\otimes w_{k_1})
      \otimes\dots\otimes
      (u_{i_\beta}\otimes v_{j_\beta}\otimes w_{k_\beta})
      \otimes h\otimes h\otimes\dots\,\in \bfH(\beta)
      ,
$$
$$
\xi':=(u_{i'_1}\otimes v_{j'_1}\otimes w_{k'_1})
      \otimes\dots\otimes
      (u_{i'_\alpha}\otimes v_{j'_\alpha}\otimes w_{k'_\alpha})
      \otimes h\otimes h\otimes\dots\in\, \bfH(\alpha)
.
$$
Then we get the same formula (\ref{eq:spherical}),
but a set of summation changes. Now we are obliged to
write $(i_s,j_s,k_s)$ (according colours)
 on edges of the white triangle with labels $s\le \beta$
 and $(i'_\sigma,j'_\sigma,k'_\sigma)$
 on the black triangle with label $\sigma\le \gamma$.
 The products are given over non-labeled triangles.

\smallskip


{\bf\punct Super-tensor products.} The previous construction
admits an extension. Let each Hilbert  space $U$, $V$, $W$
be $\Z_2$-graded,
$$
U=U_\0\oplus U_\1,\qquad
V=V_\0\oplus V_\1,\qquad
W=W_\0\oplus W_\1
,
$$
sums are orthogonal.

Take  bases $u_i\in U$, $v_j\in V$, $W_k\in W$
compatible with  gradations (i.e., $u_i\in U_{\0}$
or $U_{\1}$). Take an {\it even vector}
$$h\in H=:U\otimes V\otimes W
,$$
i.e,
$$
h\in (U_\0\otimes V_\0\otimes W_\0)\,\oplus\,
(U_\0\otimes V_\1\otimes W_\1)\,\oplus\,
(U_\1\otimes V_\0\otimes W_\1)\,\oplus\,
(U_\1\otimes V_\1\otimes W_\0)
$$
Next, we consider the tensor product
\begin{equation}
\bfH:=(H,h)\otimes (H,h)\otimes\dots
\label{eq:super}
\end{equation}
For  $\Z_2$-graded spaces $Y=Y_0\oplus Y_1$, $Z=Z_0\oplus Z_1$,
  the operator of transposition of factors
$$
Y\otimes Z \to Z\otimes Y
$$
is given by
\begin{align*}
&y_\0\,\otimes\, z_\0\,\mapsto
 \,z_\0\,\otimes\, y_\0\,\qquad
&  y_\0\,\otimes\, z_\1\, \mapsto\, z_\1\,\otimes\, y_\0,
 \\
& y_\1\,\otimes\, z_\0\,\mapsto\, z_\0\,\otimes\, y_\1,
 \qquad
& y_\1\,\otimes\, z_\1\, \mapsto\,-z_\1\,\otimes\, y_\1.
\end{align*}
This determines the action of symmetric group
in tensor products of $\Z_2$-graded spaces.

Taking the action of $\G$ in the infinite tensor product
(\ref{eq:super}, we
get a unitary representation of $\G$.

\SS


{\bf\punct Spherical functions in super-case.}
 We get the following expression
similar to (\ref{eq:spherical}):
\begin{multline}
\Phi_h(\fra)
=\\=
\sum
(-1)^\sigma
\prod\limits_{\text{black triangles $T$}}
h_{i(T)j(T)k(T)}
\prod\limits_{\text{white triangles $S$}}
\ov h_{i(S)j(S)k(S)}
\label{eq:spherical-2}
\end{multline}
We must describe a set of summation and a choose
of signs $(-1)^\sigma$.

We write basic vectors on edges of the checkers-board according
the following rules:

\smallskip

1) $u_i$ on red edges, $v_j$ on yellow edges, $w_j$ on blue edges.

\smallskip

2) The perimeter of a triangle contains even number of odd basis
elements, i.e., 0 or 2.

\smallskip

\begin{figure}

\epsfbox{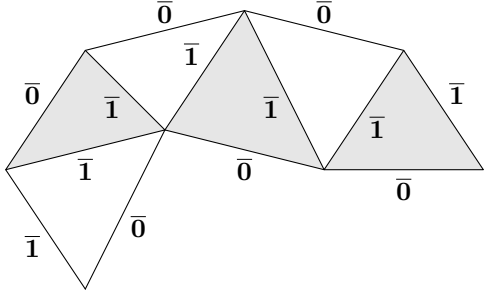}

\caption{A gallery of triangles with two odd edges}

\end{figure}

Thus, consider a triangle $T$ having 2 odd basis elements
 on its perimeter (we say 'odd edges').
 Then its neighbors through odd edges
 also have 2 odd edges. Consider their neighbors
 through odd edges etc.  In this way, we get a collection
 of closed chains of triangles. Denote by $2l_t$ their lengths
 (these numbers are even).
 Then
 $$
\sigma=\sum_t (l_t-1)
 .$$

 \section{$n$-symmetric group}

 \COUNTERS

 Here the situation is similar to the $3$-symmetric group.
All consideration given above
(checker-boards, the multiplication of double cosets,
Theorem \ref{th:io}, and constructions of representations)
 survive for  general $n\ge 2$.

 However,  there is an alternative (below {\it 'quasidual'})
  description
 of  spaces of double cosets
 $$
K(\alpha)\setminus \G^{[n]}/K(\beta)
 .
 $$
For $n>3$ the two descriptions are essentially different.

\SS


{\bf \punct Polygonal checker-boards.}
We write an element $\bfp$ of $\G_n$ as $\bfp=(p_1,\dots,p_n)$.

Consider a countable collection $A_1$, $A_2$,
\dots of oriented {\it black} $n$-gons.
On sides of each $A_k$ we write {\it anticlockwise}
 numbers $1$, $2$, \dots, $n$. Consider another
countable collection of oriented {\it white}
 $n$-gons $B_1$, $B_2$, \dots.
On sides of $B_1$, $B_2$,\dots we write {\it clockwise} numbers
$1$, $2$, \dots, $n$.

For each $i\le n$ for each $m\in\N$ we glue $i$-th side
of $A_m$ with $i$-th side of $B_{p_i(m)}$. In this
way we get an $n$-gonal complex, which is a union of countable
collection of compact closed two-dimensional surfaces.

Again, almost all components are {\it chebureks}, i.e.,
obtained by gluing of a black and white $n$-gons
along their perimeters.

\smallskip

Now we can repeat one-to-one all considerations
of Sections 2--4.

\smallskip

\begin{figure}
$$
a) \quad\epsfbox{tri.16}\qquad\qquad
 b)
\epsfbox{tri.17}
$$

\caption{A cheburek and a quasidual cheburek for $n=4$.}

a) In checker-board model a cheburek is obtained by gluing
of two quadrangles along their perimeters.
$$
\text{number of faces $=2$, number of vertices $=4$,
number of edges $=4$.}
$$

b) In the quasidual model, we glue $6$ biangles:
$$
\text{number of faces $=6$, number of vertices $=2$,
number of edges $=4$.}
$$
\label{fig:4cheburek}
\end{figure}

\begin{figure}
\epsfbox{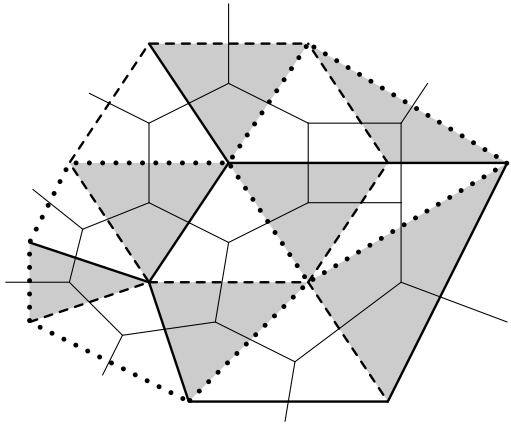}
\caption{$n=3$. A checker-board and the dual complex.}
\label{fig:dual}
\end{figure}


\begin{figure}
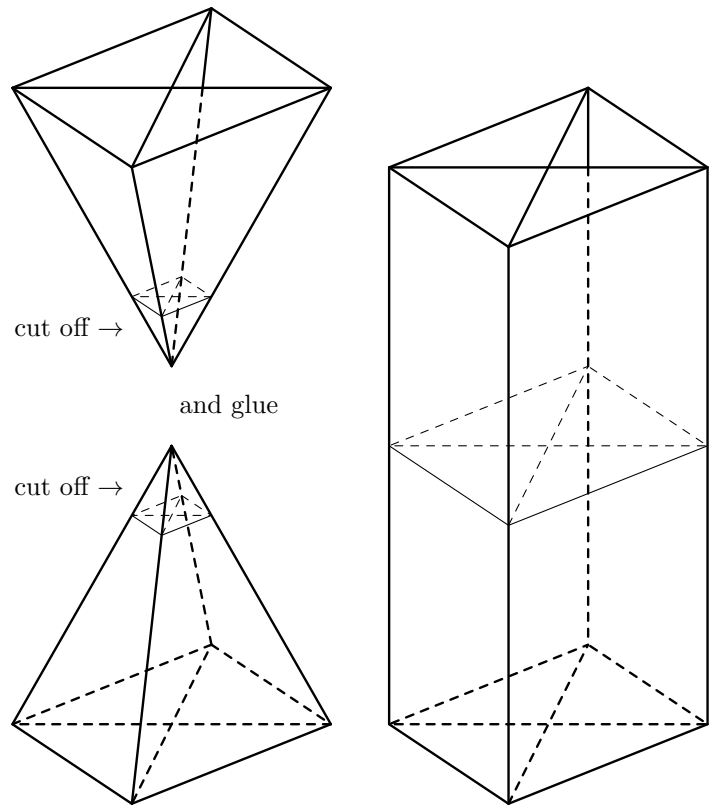

$$\epsfbox{tri.7}\qquad \epsfbox{tri.8}$$
\caption{Gluing of quasidual polygonal complexes.
 Picture at a vertex.}
\label{fig:last}
\end{figure}

{\bf \punct  The quasi-dual polygonal complex.}
 We draw a countable collection
of {\it black} vertices $a_1$, $a_2$, \dots and
a countable collection
of {\it white} vertices $b_1$, $b_2$, \dots. For
each $i\le n$ for each $m\in\N$ we
draw an edge $[a_i,b_{p_i(m)}]$ and draw the label $i$ on this edge.
In this way, we get a graph%
\footnote{Which is dual to the $1$-skeleton of the checker-board.}.

Next, for any $i$, $j\le n$, we consider the permutation
$p_{i}^{-1} p_j$. Decompose it in independent cycles.
We get  chains of the form
$$
s_1 \stackrel{p_j}{\mapsto} t_1
\stackrel{p_i^{-1}}{\mapsto}
s_2
 \stackrel{p_j}{\mapsto}
 t_2
\stackrel{p_i^{-1}}{\mapsto}s_3
\stackrel{p_j}{\mapsto}
\dots
\stackrel{p_j}{\mapsto}
t_k
\stackrel{p_i^{-1}}{\mapsto} s_1
$$

Then $s_1t_1s_2t_2\dots t_k s_1$
is a closed way on our graph. We glue a $2l$-gon to this way.
Thus
we get a polygonal complex.

\smallskip

We get $n(n-1)/2$ polygons meeting in each vertex.
Each pair of polygons has a common edge.
Therefore, {\it for $n>3$ our complex is not a
$2$-dimensional surface.} See Figure \ref{fig:4cheburek}.

\smallskip

{\sc Remark.} For $n=3$ these constructions are dual
in the usual sense, see Figure \ref{fig:dual}.
\hfill $\square$

\smallskip

To pass to double cosets $K(\alpha)\setminus \G^{[n]}/K(\beta)$,
we remove labels $>\beta$ from black vertices and
labels $>\alpha$ from white vertices.

\smallskip

It remains to describe the multiplication
$$
K(\alpha)\setminus \G^{[n]}/K(\beta)
\,\times\,
K(\beta)\setminus \G^{[n]}/K(\gamma)
\to
K(\alpha)\setminus \G^{[n]}/K(\gamma)
$$
For each $i\le \beta$ we take the link
of the black vertex
 $b_i\in K(\beta)\setminus \G^{[n]}/K(\gamma)$
and the link the white  vertex
$a_i\in K(\alpha)\setminus \G^{[n]}/K(\beta)$.
We remove neighborhoods of
 both the vertices and get collection of
half-edges in each complex. So we glue halves of edges
with the same labels together.
We must   describe two-dimensional faces of the new complex.
We glue the corresponding polygons as it is shown on
Figure \ref{fig:last}.

\bigskip

\setcounter{addendums}{1}

\section{A priori theorems about representations
of $n$-symmetric groups}

\COUNTERS

Here we use the checker-board model for double
cosets.

\smallskip


{\bf \punct Representations of products of group.
\label{ss:Aproducts}}
Let $G_1$, $G_2$ be finite groups. Then any
irreducible representation of $G_1\times G_2$ is a tensor
product $\tau_1\times\tau_2$ of  irreducible representations
of $G_1$ and $G_2$. A similar theorem for {\it unitary}
representations holds for compact groups,
semisimple groups (over $\R$ or $p$-adics), nilpotent Lie groups.
However it is not the case for unitary representations of
arbitrary groups. Apparently, a  minimal necessary condition
 (see, e.g., \cite{Dix}) is:

\smallskip

--- {\it the group $G_1$ has type I}.

\smallskip

  We can restrict an irreducible representation of $\G^{[n]}$ to
the dense subgroup
$S_\infty^{fin}\times\dots \times S_\infty^{fin}$.
But it is not a tensor product of representations
of $S_\infty^{fin}$. Moreover, consider two classes
of representation of
 $S_\infty^{fin}\times\dots \times S_\infty^{fin}$:

 \smallskip

 --- the  set of tensor products $\tau_1\otimes\dots\otimes\tau_n$;

 \smallskip

 --- the set of representations admitting a continuous extension
 to the group $\G^{[n]}$.

 \smallskip

 The intersection of these two huge classes is almost trivial.
 Precisely:

\begin{proposition}
\label{pr:Aprod}
a) Each irreducible
 representation $\tau_1\otimes\dots\otimes \tau_n$ of
$S_\infty^{fin}\times\dots \times S_\infty^{fin}$
having a  non-zero $K$-fixed vector
is one-dimensional.

\smallskip

b) Let $\rho_j$ be irreducible representations of $S_{\infty}^{fin}$.
The representation $\tau_1\otimes\dots\otimes \tau_n$
 admits a continuation to $\G^{[n]}$ if
 and only if each $\tau_j$ is continuous
 in the topology of $S_\infty$.
\end{proposition}

Proof is given in Subsection \ref{ss:Aproof}.

\smallskip


\SS

{\bf \punct Reformulations of continuity.
\label{ss:Acontinuity}}

\begin{theorem}
\label{th:continuity}
a) Let an irreducible unitary
 representation $\rho$ of
$S_\infty^{fin}\times\dots \times S_\infty^{fin}$
 have a non-zero vector invariant with respect to the
diagonal subgroup $K$. Then $\rho$ is continuous in the topology
of $\G^{[n]}$.

\smallskip

b) Let an irreducible unitary representation $\rho$ of
$S_\infty^{fin}\times\dots \times S_\infty^{fin}$ have
 a non-zero vector invariant with respect to
some   subgroup $K(\alpha)\subset K$.
 Then $\rho$ is continuous in the topology of $\G^{[n]}$.

\smallskip

c) Let $\tau$ be a unitary (reducible) representation of
$S_\infty^{fin}\times\dots \times S_\infty^{fin}$ in a Hilbert space
$H$.
Denote by $H(\alpha)$ the space of all $K(\alpha)$-fixed
vectors in $H$. Then the following conditions are equivalent:

\smallskip

--- $\tau$ is continuous in the topology of $\G^{[n]}$.

\smallskip

--- the space $\cup H(\alpha)$ is dense in $H$.
\end{theorem}

Olshanski's \cite{Olsh-symm}
 proof for $n=2$ survives for any finite $n$.

\smallskip

We formulate separately this criterion for
the group
$\G^{[1]}=S_\infty$.

\begin{theorem}
\label{th:continuity-bis}
 Let $\tau$ be a unitary (reducible) representation of
$S_\infty^{fin}$ in a Hilbert space
$H$. Then the following conditions are equivalent:

\smallskip

--- $\tau$ is continuous in the topology of $S_\infty$.

\smallskip

--- the space $\cup H(\alpha)$ is dense in $H$.

\smallskip

If $\tau$ is irreducible, then these conditions are equivalent to

\smallskip

--- the space $\cup H(\alpha)$ is non-empty.
\end{theorem}


{\bf\punct Proof of Proposition \ref{pr:Aprod}.%
\label{ss:Aproof}}

\begin{lemma}
\label{l:net}
Let $\tau_1$, $\tau_2$ be unitary (generally, reducible)
representations of a group $G$. Assume that $\tau_2$ has not
finite-dimensional $G$-invariant subspaces.
Then the representation $\tau_1\otimes\tau_2$ also
has not finite-dimensional $G$-invariant subspaces.
\end{lemma}

{\sc Proof.} Denote by $V_1$, $V_2$ the spaces of the representations
$\tau_1$, $\tau_2$. Assume that $V_1\otimes V_2$ admits
a finite-dimensional invariant subspace $W$, denote by
$\psi$ the  subrepresentation in $W$.

Then the tensor product
$\psi^*\otimes \tau_1\otimes \tau_2$ admits a $G$-invariant
vector $\eta$. Indeed,
$$
W^*\otimes (V_1\otimes V_2)\supset W^*\otimes W
.
$$
We can regard $W^*\otimes W$ as the space of operators
$W\to W$,  the identical operator
 commutes with the action
of $G$.

Next, we identify $W^*\otimes V_1\otimes V_2$
with the space of Hilbert--Schmidt operators
$$
W\otimes V_1^*\to V_2
.
$$
An invariant vector $\eta$ corresponds to an intertwining operator
$
W\otimes V_1^*\to V_2
$.
 We choose bases in the initial space
$W\otimes V_1^*$ and in the target space
$V_2$ such that the matrix of $A$ is diagonal. Equivalently,
we expand $A$ in a series
$$
A=\sum_j \lambda_j R_j,
$$
where $\lambda_j$ are the singular numbers of $A$,
and $R_j$ are partial isometries,
$$
\mathrm{im}\, R_i\bot\mathrm{im}\, R_j,\qquad
\mathrm{im}\, R_i^*\bot\mathrm{im}\, R_j^*
.$$
Then the operators $R_i$ are intertwining,
hence $\mathrm{im}\, R_j$ is an invariant subspace.
However $\rk(R_j)$ is the multiplicity of the singular
number $\lambda_j$, it is finite. Therefore
$\mathrm{im}\, R_i=0$ and $R_i=0$.
\hfill $\square$

\smallskip

{\sc Proof of Proposition \ref{pr:Aprod}.}
The group $S_{\infty}^{fin}$ has two
one-dimensional characters, the trivial one and the signature.
They are the only  irreducible
finite dimensional representations of $S_\infty^{fin}$.
Note also that the signature does not admit an extension
to the complete group $S_\infty$.

Let a representation
$\tau_1\otimes\dots\otimes\tau_n$ of
$S_\infty^{fin}\times\dots\times S_\infty^{fin}$
has a $K(\alpha)\simeq S_\infty^{fin}(\alpha)$-invariant vector.
By Lemma \ref{l:net}
 each $\tau_j$ has an $S_\infty^{fin}(\alpha)$-eigenvector
$v_j$,
$$
\tau_j(g)v_j=\chi_j(g)v_j,\qquad
g\in S_\infty^{fin}(\alpha), \qquad
\text{$\chi_j$ is a character of $S_\infty^{fin}$}
.
$$
Therefore all representations
$\wt\tau_j\simeq\chi_j\tau_j$
 have $S_\infty^{fin}(\alpha)$-invariant vectors
and are continuous in the topology of $S_\infty$.

Thus, our representation of $\G^{[n]}$ has the form
$$
(g_1,\dots,g_n)\mapsto
\Bigl(\prod\chi_j(g_j)\Bigr)
\wt\tau_1(g_1)\otimes \dots \otimes \wt\tau_n(g_n)
.
$$
Since the representations
 $$
 (g_1,\dots,g_n)\mapsto \wt\tau(g_j)
 $$
 are continuous in the topology of $\G^{[n]}$,
 it follows that the character
$\prod\chi_j(g_j)$ also is continuous. Therefore it is trivial.
\hfill $\square$

\smallskip


{\bf \punct  Analogy with $p$-adic groups.}
 The criterion  of continuity
 (Theorem \ref{th:continuity})
is an imitation of admissibility for $p$-adic groups.
For definiteness, consider $G=\GL(n,\Q_p)$. Denote by
$\GL^\alpha$ the subgroup consisting of matrices
$1+p^\alpha T$, where $T$ has integer coefficients.
Then a representation $\rho$ of $G$ in a space $H$ is {\it admissible}
if each vector is fixed by some subgroup $\GL^\alpha$.
A unitary representation of $G$ in a Hilbert space $H$
is admissible, more precisely it is admissible on a dense
subspace.

\smallskip

The reason of this parallel, see \cite{Ner-book}, Proposition
VIII.1.3.

\smallskip



{\bf \punct Correspondence between representations
of the group $\G^{[n]}$ and representations
of the category $\bbS^{[n]}$.}
In Section 3 we obtained the canonical map
$$
\Bigl\{\text{Unitary representations of $\G^{[n]}$}
\Bigr\}
\to
\Bigl\{
\text{$*$-representations of $\bbS^{[n]}$}
\Bigr\}
\,.
$$

\begin{theorem}
\label{th:bijection1}
This map is a bijection.
\end{theorem}

Below we prove a stronger version of this
statement.


\smallskip

{\bf\punct The completion of the category $\bbS^{[n]}$.}
Let us define a 'new' category $\ov \bbS^{[n]}$ obtained
 by an adding
of an  infinite object.
Its objects are $0$, $1$, $2$, \dots, $\infty$.

A morphism $\beta\to\alpha$ is the following collection
of data:

\smallskip

--- a finite or countable collection of checker-boards
colored in the usual way;

\smallskip

--- almost all components are chebureks;

\smallskip

--- we  fix an injective map from the set
 $\{1,2,\dots,\beta\}$ to the set of black faces
 and an injective map  $\{1,2,\dots,\alpha\}$
  to the set of white faces;

\smallskip

--- there is no empty chebureks.

\smallskip

The product of morphisms is given in the usual way.

\smallskip

We use the common notation $\Mor(\beta,\alpha)$ for sets
of morphisms $\beta\to\alpha$, $\End(\alpha)$ for semigroups
of endomorphism, and $\Aut(\alpha)$ for group
of automorphisms.

Evidently,
$$\Aut(\alpha)=S_\alpha\times\dots\times S_\alpha,
\,\,\text{where $\alpha<\infty$},
 \qquad \Aut(\infty)=\G^{[n]}
 .
 $$

\smallskip

At the end of this subsection we define the natural topology
on the semigroup $\End(\infty)$ (representations
must be continuous%
\footnote{With respect to the weak  topology
on the space of operators in a Hilbert space})
 and discuss some properties
of the category $\ov\bbS^{[n]}$. Now we formulate
the main statement of the subsection.

\begin{theorem}
\label{th:Agk2}
 Any $*$-representation of the category $\bbS^{[n]}$
admits a unique continuous extension to
the $*$-representation of the category $\ov\bbS^{[n]}$.
\end{theorem}

 We need the following distinguished morphisms
$$
\text{$\lambda_{\beta,\alpha}:\beta\to\alpha$,\qquad
$\mu_{\alpha,\beta}:\alpha\to\beta$,\qquad
$\theta^\alpha_\beta:\alpha\to\alpha$
}
$$
defined for all $\alpha>\beta$.

The morphism  $\lambda_{\beta,\alpha}:\beta\to\alpha$
is defined as a
 disjoint sum of  chebureks with the following labels
 on black/white sides:
\begin{gather*}
\text{$(j/j)$, where $j\le\beta$,}
\\
\text{$(p/\text{empty})$,
where $\beta<p\le\alpha$}.
\end{gather*}
We set
$$
\mu_{\alpha,\beta}:=
\lambda_{\beta,\alpha}^\square
,\qquad
\theta^\alpha_\beta:=\lambda_{\beta,\alpha}\mu_{\alpha,\beta}
$$
Thus $\mu_{\alpha,\beta}$ is determined by the following
collection of chebureks:
\begin{gather*}
\text{$(j/j)$, where $j\le\beta$,}
\\
\text{ $(\text{empty}/p)$,
where $\beta<p\le\alpha$}.
\end{gather*}
and $\theta^\alpha_\beta$ by
\begin{gather}
\label{A.1}
\text{$(j/j)$, where $j\le\beta$,}
\\
\label{A.2}
\text{$(p/\text{empty})$, $(\text{empty}/p)$,
where $\beta<p\le\alpha$}.
\end{gather}

We need some simple fact concerning these morphisms.

First,
\begin{gather}
\label{A.3}
(\theta^\alpha_\beta)^2=\theta^\alpha_\beta, \qquad
\theta^\alpha_\beta=(\theta^\alpha_\beta)^\square
\\
\label{A.4}
\lambda_{\beta,\alpha}^\square \lambda_{\beta,\alpha}=1
\end{gather}

Second, the map
\begin{equation}
\label{A.5}
\fra\mapsto \lambda_{\beta,\alpha} \fra\mu_{\alpha,\beta}
\end{equation}
is an injective homomorphism of the semigroup $\End(\beta)$
to the semigroup $\End(\alpha)$. In fact, we add
chebureks (\ref{A.2}) to the checker-board $\fra$.

Third, we need the map
$
\End(\alpha)\to\End(\beta)
$
given by
$$
\frb\mapsto \mu_{\alpha,\beta}\frb \lambda_{\beta,\alpha}
$$
This is equivalent to removing of labels $>\beta$
from black and white triangles.

Now we are ready to define the topology on $\End(\infty)$.
We say that a sequence $\frc_j\in\End(\infty)$ converges
to $\frc$ if for each $\beta<\infty$
$$
\mu_{\infty,\beta}\frc_j \lambda_{\beta,\infty}
=
\mu_{\infty,\beta}\frc \lambda_{\beta,\infty},
\qquad
\text{for sufficiently large $j$.}
$$

{\sc Example.} A weak zero $h_jK(\alpha)$,
see Proposition \ref{ss:weak}, converges
to $\theta^\infty_\alpha$. \hfill $\square$

\begin{proposition}
a) The multiplication in $\End(\infty)$
is separately continuous.

\smallskip

b) The group $\Aut(\infty)$ is dense in the semigroup
 $\End(\infty)$.
\end{proposition}

{\sc Proof of b.}
Let $\frc\in\End(\infty)$
have non-labeled white and non-labeled
black faces. For each $\beta$ let us make
$\frd_\beta\in \End(\infty)$, where all
black faces are labeled, in the following way.

We preserve all white labels on $\frc$,
preserve  all black label $\le\beta$, remove black labels
with numbers $>\beta$, and enumerate
all initially non-labeled black faces
and all faces liberated from labels in arbitrary way.
Then $\frd_\beta$ converges to $\frc$.

Next, for each $\frd_\beta$ we construct a sequence
$\fre_{\beta\alpha}$ by repeating the same actions
for white faces. Then $\fre_{\beta\alpha}$
converges to $\frd_\beta$.

On the other hand, $\fre_{\beta\alpha}\in\End(\infty)$.
\hfill $\square$


{\bf\punct A formal proof of Theorem \ref{th:Agk2}.%
\label{ss:formal}}
It is a special case of
 Theorem VIII.1.10 from \cite{Ner-book}.
All necessary structures were defined in the previous subsection.
We must only verify:

\begin{lemma}
For any $*$-representation $\ov\rho$ of the category
$\bbS^{[n]}$ for any morphism
$\fra$,
$$
\|\ov\rho(\fra)\|\le 1.
$$
\end{lemma}

{\sc Proof.} By (\ref{A.3})--(\ref{A.4}),
$\ov\rho(\theta^\alpha_\beta)$ is an orthogonal projection
and $\ov\rho(\lambda_{\beta,\alpha})$ is an isometric embedding.
For $\frc\in\Aut(\gamma)$, the operator
$\ov\rho(\frc)$ is unitary.

 We can represent any
$\fra \in\Mor (\beta,\alpha)$ as
$$
\fra=\mu_{\gamma,\alpha} \frc \nu_{\beta,\gamma},
\qquad \frc\in \Aut(\gamma)
$$
for some $\gamma$. We simply write labels on all empty
triangles. This completes the proof.
\hfill $\square\square$

\smallskip



\smallskip

{\bf \punct Reconstruction of a representation of the
group $\G^{[n]}$ from a representation of the category.%
\label{ss:reconstruction}}
We will not repeat the proof of
Theorem VIII.1.10 from \cite{Ner-book}.
Here we briefly describe  the limit pass
necessary for Theorem \ref{th:Agk2}.

Let $\rho$ be a $*$-representation  of $\bbS^{[n]}$.
Let $H(\alpha)$ be the spaces of representation.
For $\beta<\alpha$, the operator $\rho(\lambda_{\beta,\alpha})$
is an isometric embedding $H(\beta)\to H(\alpha)$.
The semigroup $\End(\beta)$ acts in the both spaces,
see (\ref{A.5}); the operator $\rho(\lambda_{\beta,\alpha})$
intertwines these actions.

Therefore, we can assume that $H(\beta)\subset H(\alpha)$
and take the inductive limit $H(\infty)$ of the chain
$$
\dots \to H(\beta)\to H(\beta+1)\to \dots
$$
of Hilbert spaces (i.e., we take the completion of
$\cup_\beta H(\beta)$).
The inductive limit
$\cup_{\beta<\infty} \End(\beta)$
 of the chain of semigroups
$$
\dots\to \End(\beta)\to \End(\beta+1)\to\dots
$$
acts in $H(\infty)$.

However, the semigroup $\cup_{\beta<\infty} \End(\beta)$
is a small subset in $\End(\infty)$
($\cup_{\beta<\infty} \End(\beta)$ is countable and   contains no
invertible elements at all).
 To define  $\ov\rho(\frc)$ for an arbitrary
$\frc\in\End(\infty)$, we write
$\ov\rho(\frc)$ as a weak limit
\begin{equation}
\ov\rho(\frc):=
\lim_{\beta\to\infty}
\ov\rho(\theta^\infty_\beta\frc \theta^\infty_\beta )
.
\label{eq:chain1}
\end{equation}


{\bf\punct The composition of functors.}
Thus we get a chain of 3 functors:
\begin{multline}
\Bigl\{\begin{matrix}\text{Unitary}\\
\text{representations of $\G^{[n]}$}
\end{matrix}
\Bigr\}
\stackrel{T_1}
{\longrightarrow}
\Bigl\{
\begin{matrix}
\text{$*$-Representations}\\
\text{of $\bbS^{[n]}$}
\end{matrix}\Bigr\}
\stackrel{T_2}\longrightarrow
\\
\longrightarrow
\Bigl\{
\begin{matrix}
\text{$*$-Representations}\\
\text{of $\ov\bbS^{[n]}$}
\end{matrix}\Bigr\}
\stackrel{T_3}
\longrightarrow
\Bigl\{\begin{matrix}\text{Unitary}\\
\text{representations of $\G^{[n]}$}
\end{matrix}
\Bigr\}
\label{eq:functors}
\end{multline}
The last map is the restriction of a representation of
$\End(\infty)$ to $\Aut(\infty)\simeq \G^{[n]}$.

We must show that this chain is closed,

\begin{proposition}
The composition
$T_3 T_2 T_1$ is the identical map.
\label{pr:id}
\end{proposition}

{\sc Proof.}
Let $\rho$ be a unitary representation of $\G^{[n]}$.
Let $\frb\in\End(\beta)$, where $\beta<\infty$.
In Section 3, we  defined the operator
$\ov\rho(\frb)$ as follows. Take $\bfg_\beta\in \G^{[n]}$
such that
\begin{equation}
\frb=\theta^\infty_\beta \bfg_\beta \theta^\infty_\beta
.
\label{eq:chain2}
\end{equation}
Then
\begin{equation}
\ov\rho(\frb):=P(\beta)\rho(\bfg_\beta) P(\beta)
.
\label{eq:chain3}
\end{equation}
Take $\bfg\in\Aut(\infty)\simeq\G^{[n]}$.
Keeping in the mind
(\ref{eq:chain1})--(\ref{eq:chain3}), we get
\begin{equation}
\ov\rho(\bfg)=\lim_{\beta\to\infty}
 P(\beta)\rho(\bfg)P(\beta)
\label{eq:chain4}
\end{equation}
(the limit is a weak limit).
By Theorem \ref{th:continuity}, $P(\beta)$
strongly tends to the identical operator.
Therefore the weak limit (\ref{eq:chain4}) is $\rho(\bfg)$.

 Thus,
 $
\ov \rho(\bfg)= \rho(\bfg)$ for $\bfg\in\G^{[n]}$.
 \hfill $\square$

 \smallskip

Next, we prove a stronger form of Proposition
\ref{pr:id}.

\begin{theorem}
\label{th:id}
All the maps $T_1$, $T_2$, $T_3$ are bijections.
The composition $T_3 T_2 T_1$,
 $T_2 T_1 T_3$, $T_1 T_3 T_2$
are identical.
\end{theorem}

{\sc Proof.} {\it The injectivity of $T_1$.}
This follows from the identity
$$\lim_{\beta\to\infty}
 P(\beta)\rho(\bfg)P(\beta)
 =
\rho(\bfg)
,$$
see the previous proof.

\smallskip

{\it The bijectivity of $T_2$} is evident,
the map $T_2^{-1}$ is the restriction of a representation
of $\ov\bbS^{[n]}$ to $\bbS^{[n]}$.

\smallskip

{\it The injectivity of $T_3$.} The group
$\Aut(\infty)$ is dense in the semigroup $\End(\infty)$.
Therefore a representation of $\Aut(\infty)$ remembers
a representation of $\End(\infty)$. On the other hand, all
the semigroups $\End(\alpha)$ can be regarded as subsemigroups
in $\End(\infty)$.
\hfill
$\square$


\smallskip

{\bf\punct Uniqueness of a spherical vector.}

\begin{theorem}
Let $\rho$ be a unitary irreducible representation
of $\G^{[n]}$, let $H(0)$ be the space of $K$-fixed vectors.
Then
$$
\dim H(0)\le 1.
$$
\end{theorem}

{\sc Proof.} Our functors (\ref{eq:functors})
 send direct sums
to direct sums and irreducible representations to irreducible%
\footnote{For categories, see definitions in
\cite{Ner-book}, II.8.}
representations. In particular,
the representation of
$\End(0)\simeq K(0)\setminus \G^{[n]}/K(0)$
must be irreducible (see, e.g., \cite{Ner-book},
Lemma II.8.1.). The semigroup
$\simeq K(0)\setminus \G^{[n]}/K(0)$
is a commutative semigroup with involution.
Its irreducible $*$-representations are one-dimensional.
\hfill
$\square$

\smallskip

If $\dim H(0)=1$, we say that the representation is
{\it $K$-spherical}.


\smallskip


{\bf \punct Spherical characters.}
For $\alpha=0$, $1$, \dots, $\infty$,
denote by $\Delta(\alpha)$ the center of $\End(\infty)$.

\begin{proposition}
All elements of $\Delta(\alpha)$ have the form:
the collection of labeled chebureks
$$
(j/j),\quad \text{where $j\le\alpha$},
$$
and a finite disjoint union of compact non-labeled
surfaces.
\end{proposition}

The statement is obvious.

\smallskip

Thus, all semigroups $\Delta(\alpha)$ are isomorphic
to
$$
\Delta:=\Delta(0)=K(0)\setminus \G^{n]}/K(0)
.
$$

For $\alpha>\beta$, the isomorphism $\Delta(\alpha)\to\Delta(\beta)$
can be written as follows:
\begin{equation}
\frz\mapsto\mu_{\alpha,\beta}  \frz \lambda_{\beta,\alpha}
\label{eq:center-map}
.\end{equation}

The center of $\End(\infty)$ acts in  irreducible representations
by scalar operator.

\begin{observation}
For any irreducible representation $\rho$ of $\G^{[n]}$,
we get a homomorphism $\chi_\rho$
 of $\Delta$ to the multiplicative semigroup
of $\C$.
\end{observation}

We say that $\chi_\rho$ is the {\it spherical character} of
the representation $\rho$.

\begin{observation}
For any irreducible spherical representation of $\G^{[n]}$
the spherical character coincide with spherical function.
\end{observation}

\begin{lemma}
For any irreducible representation
$\rho$ of
$\G^{[n]}$
for each $\alpha$ the semigroup $\Delta(\alpha)$ acts
in $H(\alpha)$ as $\chi_\rho$.
\end{lemma}

{\sc Proof.} We apply the map (\ref{eq:center-map}).
\hfill$\square$

\smallskip


{\bf \punct Self-similarity and spherical characters.%
\label{ss:Aself}}
 Consider the subgroup $\G^{[n]}(\alpha)\subset \G^{[n]}$
that fixes points of the initial segment $1$, \dots, $\alpha$
 in all $n$ copies of $\N$. Then
$$
\G^{[n]}(\alpha)
\simeq
\G^{[n]}
.
$$
Restricting a representation of $\G^{[n]}$
to $\G^{[n]}(\alpha)$, we get a collection
of (reducible) representations of $\G^{[n]}$.

\smallskip

Evidently, the embedding
$$\G^{[n]}(\alpha)
\to\G^{[n]}$$
induces an isomorphism of centers of two copies of
$\End(\infty)$. This implies the following statement:

\begin{theorem}
 Let $\rho$ be an irreducible unitary representation
of $\G^{[n]}$. Consider the restriction
$$
\rho\Bigr|_{\G^{[n]}(\alpha)}
$$
and a spherical subrepresentation $\tau$ of the restriction%
\footnote{For sufficiently large $\alpha$ such  subrepresentations
exist.}.
Then the spherical function of $\tau$ coincides
 with the spherical character of $\rho$.
\end{theorem}

{\sc Example.}  A representation of $\G^{[n]}$ is continuous
in the topology of $S_\infty\times\dots\times S_\infty$
(see Subsection 6.1) if and only if its spherical character $=1$
(this is a rephrasing of Proposition \ref{pr:Aprod}).
\hfill
$\square$

\smallskip

Let $\chi$ be a spherical character. For each
connected non-labeled checker-board $\frc$ we have a number
$\chi(\frc)$. For a disjoint sum $\frc_1\cup\dots\cup \frc_m$
we have
$$
\chi(\frc_1\cup\dots\cup \frc_m)=\prod_j\chi(\frc_j)
.
$$
Thus a spherical character is determined by numbers
$\chi(\frc)$, where $\frc$ ranges in set of connected checker-boards.

However, numbers $\chi(\frc)$ are not arbitrary.

\smallskip


{\bf \punct Status of problem of classification of
unitary representations of the group $\G^{[n]}$.%
\label{ss:status}}
It is not clear, is this problem reasonable or not
(there are different arguments 'pro and con').

However, the problem looks like two-step:

\smallskip

a) to describe all possible spherical characters%
\footnote{The construction of Section 4 does not exhaust
all spherical characters. Even for $\G^{[2]}$ there is an additional
parameter.};

\smallskip

b) to classify all irreducible representations
with a given spherical character.

\smallskip

There is also (apparently, simple) problem:

\smallskip

c) to prove that $\G^{[n]}$ is a type I group.

\smallskip

 For
the bisymmetric group $\G^{[2]}$, the description of
spherical characters
is equivalent to the famous Thoma Theorem (1964),
 see Theorem \ref{th:Thoma}.
The problem b) was solved by Olshanski \cite{Olsh-symm}
and Okounkov \cite{Oko}.

\smallskip


\begin{figure}
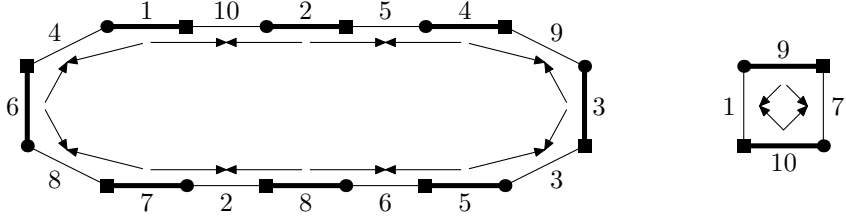

\epsfbox{tri.10}
\caption{Example. A pair of substitutions}
\label{fig:one-dim}
$$
\begin{pmatrix}
1 & 2 &  4 & 3 & 5 & 8 & 7 & 6& 9& 10\\
4 & 10 & 5 & 9 & 3 & 6 & 2 & 8& 1 &7
\end{pmatrix},
\begin{pmatrix}
1 & 2 &  4 &  3 &  5 & 8 & 7 & 6 & 9 & 10\\
10 & 5 &  9 & 3 & 6 & 2 & 8 & 4 & 7 & 1
\end{pmatrix}
$$
Notation for colors:

\epsfbox{tri.11}
\end{figure}

\begin{figure}
$$\epsfbox{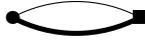}$$

\caption{The analog of a cheburek}
\label{fig:cheburek-bis}

\end{figure}

\begin{figure}

\caption{The production of a chip from an one-dimensional
simplicial complex.}
\label{fig:chips}

We take an element of $K(4)\setminus\G^{[2]}/K(3)$:
$$
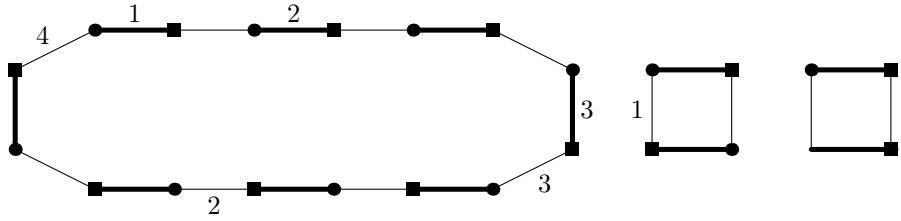
\epsfbox{tri.12}
$$
Cut labeled edges:
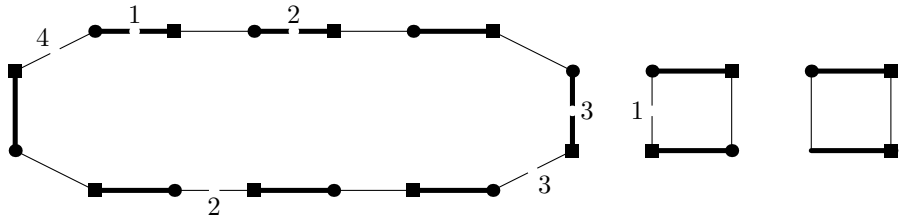
$$
\epsfbox{tri.13}
$$
We get a collection of arcs and draw this collections as a chip:
$$
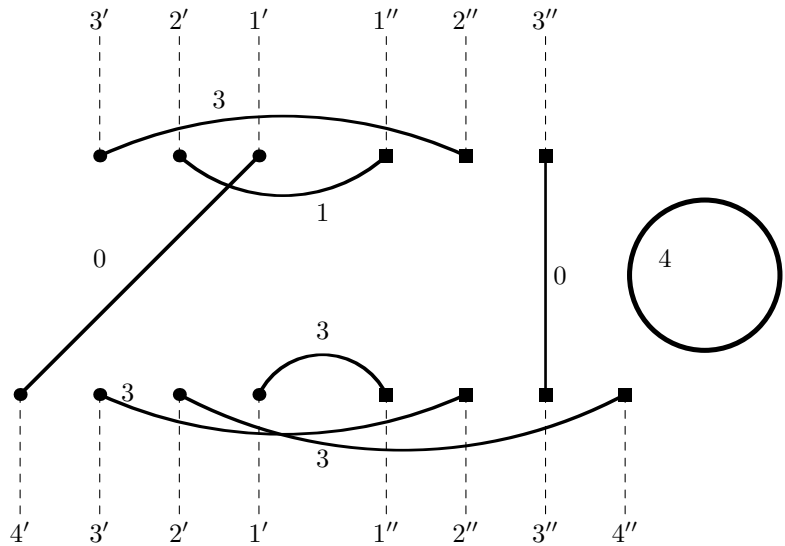
\epsfbox{tri.15}
$$
\end{figure}

\bigskip

\newpage

{\bf\large A. Addendum. Bisymmetric group}

\bigskip

{\bf \Apunct  Olshanski chips.%
\label{ss:Achips}} Let $n=2$.
Consider an element $\bfg=(g_{red},g_{blue})\in \G^{[2]}$.
Take a countable collection of black segments $A_j$
with red and blue ends. Take a countable collection
of white segments with red and blue ends.
If $g_{red} i=j$, then we glue the red end of $A_i$
with the red end of $B_j$. If $g_{blue}k=l$, then
we glue the blue end of $A_k$ with the blue end of $B_l$.
In this way, we get an one-dimensional complex.
One-dimensional strata of the complex (segments)
 are colored in black and
white, zero-dimensional strata in red and blue.
 There are labels $\in \N$
on black segments and labels $\in\N$ on white segments.

The pair $\bfg=(g_{red},g_{blue})$ can be easily reconstructed
 from these data, see Figure \ref{fig:one-dim}.

 To realize a double coset
 $$
\in K(\alpha)\setminus\G^{[2]}/K(\beta)
 ,$$
we remove labels $>\beta$ from black segments and labels
$>\alpha$ from white segments. The analog of an empty cheburek
is a circle consisting of two non-labeled segments.
It is reasonable to throw out all empty chebureks
(see Figure \ref{fig:cheburek-bis}).

\smallskip

Thus we get a finite disjoint union of closed chains of segments,
black and white segments  interlace, red and blue ends also interlace.
There are pairwise different labels $1$, \dots, $\beta$
on some black segments and pairwise different labels
$1$, \dots, $\alpha$ on some white segments.

\smallskip

Next, it is reasonable to do the following operation:

\smallskip

1. We cut  all labeled segments at midpoints
  and  assign the corresponding labels to the ends.

\smallskip

2.  After this we get a finite collection of cycles having
  even lengths $2 l_p$ and a finite collection of non-closed chains.
  If both ends of a chain are white (resp. black), then
  the chain contains an odd number $2m_q+1$ of complete segments.
  If one end is black and another is white, then the chain
  contains an even  number $2k_r$ of complete segments.
  See Figure \ref{fig:chips}.

\smallskip

Now an element of $K(\alpha)\setminus\G^{[2]}/K(\beta)$
corresponds to diagram (a {\it chip}) of the following form
(see Figure \ref{fig:chips}):

\smallskip

a) There are $\beta+\beta$ {\it entries}. On Figure
we show them by $\beta$ labeled
fat points and  $\beta$ labeled square points.
There are $\alpha+\alpha$ {\it exits}. On Figure \ref{fig:chips}
we show them by numbered $\alpha$
fat points and numbered $\alpha$ square points.

\smallskip

b) There are arcs of 3 types (horizontal, vertical, and cyclic).
 For each arc we assign a non-negative
integer (length). The rules of designing are:

\smallskip

--- a horizontal arc connects a square point
and a fat point in the entrance
or a square and a fat point in the exit. It has odd length.

\smallskip

--- a vertical arc connects a fat point in the entrance with
a fat point in the exit or a square point in the entrance with
a square point in the exit. It has even length.

\smallskip

--- cycles have even length.

\smallskip

Having two chips $\fra:\gamma\to\beta$ and $\frb:\beta\to\alpha$.
Then we connect exits of $\fra$ with corresponding
entrances of $\frb$ and evaluate lengths of  composite
arcs.

\smallskip

\begin{figure}
$$
\epsfbox{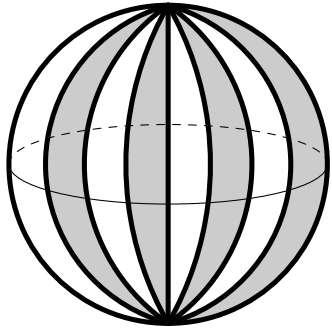}
$$
\caption{A bianglular checker-board. We omit colors on edges.}
\label{fig:biangle}

\end{figure}

  {\sc Remark.} Certainly, we can apply the general
  construction of checker-board for $\G^{[2]}$.
   In this case 2-dimensional complexes
  are spheres divided into biangles, see Figure \ref{fig:biangle}.

\smallskip

{\bf \Apunct Characters of $S_\infty$ and the bisymmetric group.%
\label{ss:AThoma}}
In 1964 E.~Thoma obtained the classification of all characters
of $S_\infty^{fin}$. According his definition,  a {\it character} is
an extreme point of the set of all functions $F$ on
$S_\infty^{fin}$ satisfying the conditions:

\smallskip

--- $F$ is central, i.e., $F(hgh^{-1})=F(g)$;

\smallskip

--- $F$ is positive definite;

\smallskip

--- $F(1)=1$.

\begin{Atheorem}   {\bf(Thoma, \cite{Tho})}
\label{th:Thoma}
All characters
of $S_{\infty}^{fin}$
have the form

\begin{equation}
\chi(g)=
\prod_{k\ge 2}
\Bigl(\sum_j \alpha_j^k-\sum_j (-\beta_j)^{k}
\Bigr)^{r_k(g)}
,
\tag{A.1}
\end{equation}
where $r_k(g)$ is the number of cycles of  length
$k$ in $g$, and the parameters $\alpha$, $\beta$ satisfy
$$
\alpha_1\ge\alpha_2\ge\dots \ge 0,\qquad
\beta_1\ge\beta_2\ge\dots \ge 0,
\qquad
\sum \alpha_j +\sum \beta_j\le 1
.
$$
\end{Atheorem}

Having a positive definite function $\chi$
on the group $S_\infty^{fin}$, we define a Hilbert space
and a unitary representation of the group in the usual way, see,
e.~g., \cite{Dix}. This representation
 is a $II_1$-factor representation
 (in particular, it is strongly reducible).

\begin{Atheorem} (\cite{Olsh-symm}).
a) There is a canonical one-to-one correspondence between
Thoma characters and spherical representations of the
bisymmetric group.

\smallskip

b) Moreover, for a character $\chi(g)$ the corresponding spherical
function is given by
$$
\Phi(g_1,g_2)=\chi(g_1 g_2^{-1})
.
$$
\end{Atheorem}

Formula (A.1) is simpler than our formulas (\ref{eq:spherical})
and (\ref{eq:spherical-2}). In fact a vector
$$
h\in U\otimes V, \qquad \text{where $U$, $V$ are Hilbert spaces,}
$$
can be reduced to the form
$$
h=\sum \alpha_j^{1/2}\, u_j\otimes v_j
$$
by a change of bases in $U$ and $V$. For a triple tensor
product such reduction is impossible.

\smallskip


{\bf  \Apunct Symmetric pairs.} Olshanski
researches in infinite-dimensional
groups (see \cite{Ols-howe}, \cite{Olsh-symm})
 partially are an imitation of symmetric pairs
and the theory of spherical functions. In Neretin \cite{Ner-book},
VIII.5 it was observed that some tricks related to infinite
dimensional symmetric pairs can be applied in a wider
 generality.


{\tt Math.Dept., University of Vienna,

 Nordbergstrasse, 15,
Vienna, Austria

\&

Institute for Theoretical and Experimental Physics,

Bolshaya Cheremushkinskaya, 25, Moscow 117259,
Russia

\&

Mech.Math. Dept., Moscow State University,
Vorob'evy Gory, Moscow


e-mail: neretin(at) mccme.ru

URL:www.mat.univie.ac.at/$\sim$neretin

wwwth.itep.ru/$\sim$neretin
}

\end{document}